\documentclass[11pt,leqno]{amsart}

\usepackage{amsthm}
\usepackage{amsmath}
\usepackage{amssymb}
\usepackage[all,cmtip]{xy}
\usepackage[T1]{fontenc}
\usepackage{mathrsfs}
\usepackage{enumerate}

\usepackage{color}

\theoremstyle{plain}
\newtheorem{thm}{Theorem}
\newtheorem{lem}{Lemma}
\newtheorem{prop}{Proposition}
\newtheorem{cor}{Corollary}
\theoremstyle{definition}
\newtheorem{defn}{Definition}

\newtheorem{exmp}{Example}

\DeclareMathOperator{\ob}{ob}
\DeclareMathOperator{\mor}{mor}

\DeclareMathOperator{\Endom}{End}

\title{Noncrossed Product Matrix Subrings and Ideals of Graded Rings}
\date{}

\author{Johan \"{O}inert and Patrik Lundstr\"{o}m}

\address{Johan \"{O}inert,
Centre for Mathematical Sciences,
Lund University,
P.O. Box 118,
SE-22100 Lund,
Sweden}
\email{Johan.Oinert@math.lth.se}

\address{Patrik Lundstr\"{o}m,
University West,
Department of Engineering Science,
SE-46186 Trollh\"{a}ttan,
Sweden}
\email{Patrik.Lundstrom@hv.se}

\begin{document}

\maketitle


\begin{abstract}
We show that if a groupoid graded ring has
a certain nonzero ideal property and
the principal component of the ring is commutative,
then the intersection of a
nonzero twosided ideal of the ring with the commutant
of the principal component of the ring is nonzero.
Furthermore, we show that for a skew groupoid ring with commutative principal component, the principal component is maximal commutative if and only if it is intersected nontrivially by each nonzero ideal of the skew groupoid ring.
We also determine the center of strongly groupoid
graded rings in terms of an action on the ring
induced by the grading.
In the end of the article, we show that, given
a finite groupoid $G$, which has a nonidentity
morphism, there is a ring, strongly graded
by $G$, which is not a crossed product over $G$.
\end{abstract}

\section{Introduction}

Let $R$ be a ring.
By this we always mean that $R$ is an additive group
equipped with a multiplication which is associative and unital.
The identity element of $R$ is denoted $1_R$ and
is always assumed to be nonzero.
We always assume that ring homomorphisms respect the multiplicative identities.
The set of ring endomorphisms of $R$ is denoted $\Endom(R)$ and
the center of $R$ is denoted $Z(R)$.
By the commutant
of a subset $S$ of a ring $R$, denoted $C_R(S)$, we mean the set of elements
of $R$ that commute with each element of $S$.

Suppose that $R_1$ is a subring of $R$ i.e. that there is
an injective ring homomorphism $R_1 \rightarrow R$.
Recall that if $R_1$ is commutative, then it
is called a maximal commutative subring of $R$
if it coincides with its commutant in $R$.
A lot of work has been devoted to investigating the
connection between on the one hand maximal
commutativity of $R_1$ in $R$ and on the other hand
nonemptyness of intersections of $R_1$ with
nonzero twosided ideals of $R$
(see \cite{coh}, \cite{fis}, \cite{for78}, \cite{lor},
\cite{lor79}, \cite{lor80},
\cite{mon78}
and \cite{pas77}).
Recently (see \cite{oin06}, \cite{oin07}, \cite{oin08},
\cite{oin09} and \cite{oin10})
such a connection was established for the commutant $R_1$
of the identity component of strongly group graded rings
and group crossed products (see Theorem \ref{TheoremOne} and Theorem \ref{TheoremTwo} below).
Let $G$ be a group with identity element $e$.
Recall that a ring $R$ is graded by the group $G$
if there is a set of additive subgroups,
$R_s$, $s \in G$, of $R$ such that $R =
\bigoplus_{s \in G} R_s$ and $R_s R_t \subseteq
R_{st}$, $s,t \in G$; if $R_s R_t = R_{st}$,
$s,t \in G$, then $R$ is called strongly graded.
The subring $R_e$ of $R$ is called the identity component
of $R$.
The following result appears in \cite{oin10}.

\begin{thm}\label{TheoremOne}
If a strongly group graded ring has commutative
identity component, then the intersection of a
nonzero twosided ideal of the ring with the commutant
of the principal component in the ring is nonzero.
\end{thm}

Recall that crossed products are defined by
first specifying a crossed system i.e. a quadruple
$\{ A,G,\sigma,\alpha \}$ where $A$ is a ring, $G$ is a group and
$\sigma : G \rightarrow \Endom(A)$ and $\alpha : G \times G \rightarrow
A$ are maps satisfying the following four conditions:
\begin{equation}\label{id}
\sigma_e = {\rm id}_A
\end{equation}
\begin{equation}\label{identity}
\alpha(s,e) = \alpha(e,s) = 1_A
\end{equation}
\begin{equation}\label{associative}
\alpha(s,t) \alpha(st,r) = \sigma_s(\alpha(t,r)) \alpha(s,tr)
\end{equation}
\begin{equation}\label{algebra}
\sigma_s(\sigma_t(a)) \alpha(s,t) = \alpha(s,t) \sigma_{st}(a)
\end{equation}
for all $s,t,r \in G$ and all $a \in A$.
The crossed product, denoted
$A \rtimes_{\alpha}^{\sigma} G$, associated to this quadruple,
is the collection of formal sums
$\sum_{s \in G} a_s u_s$, where $a_s \in A$, $s \in G$,
are chosen so that all but finitely many of them are zero.
By abuse of notation we write $u_s$ instead of $1u_s$
for all $s \in G$.
The addition on $A \rtimes_{\alpha}^{\sigma} G$ is defined pointwise
\begin{equation}\label{addition}
\sum_{s \in G} a_s u_s + \sum_{s \in G} b_s u_s =
\sum_{s \in G} (a_s + b_s)u_s
\end{equation}
and the multiplication on $A \rtimes_{\alpha}^{\sigma} G$ is defined
by the bilinear extension of the relation
\begin{equation}\label{multiplication}
(a_s u_s)(b_t u_t) = a_s \sigma_s(b_t) \alpha(s,t) u_{st}
\end{equation}
for all $s,t \in G$ and all $a_s,b_t \in A$.
By (\ref{id}) and (\ref{identity}) $u_e$ is a multiplicative identity of
$A \rtimes_{\alpha}^{\sigma} G$ and by (\ref{associative})
the multiplication on $A \rtimes_{\alpha}^{\sigma} G$
is associative. There is also an $A$-bimodule structure on $A \rtimes_{\alpha}^{\sigma} G$
defined by the linear extension of the relations
$a(b u_s) = (ab) u_s$ and
$(a u_s)b = (a \sigma_s(b)) u_s$ for
all $a,b \in A$ and all $s,t \in G$,
which, by (\ref{algebra}), makes
$A \rtimes_{\alpha}^{\sigma} G$ an $A$-algebra.
Note that $A \rtimes_{\alpha}^{\sigma} G$
is a group graded ring with the grading
$(A \rtimes_{\alpha}^{\sigma} G)_s = Au_s$, $s \in G$;
it is clear that this makes $A \rtimes_{\alpha}^{\sigma} G$
a strongly graded ring if and only if each
$\alpha(s,t)$, $s,t \in G$, has a left inverse in $A$.
In \cite{oin09}, the following result was shown.

\begin{thm}\label{TheoremTwo}
If $A \rtimes_{\alpha}^{\sigma} G$ is a crossed product
with $A$ commutative, all $\sigma_s$, $s \in G$,
are ring automorphisms and none of the $\alpha(s,s^{-1})$, $s \in G$,
are zero divisors in $A$, then every
intersection of a nonzero twosided ideal of $A \rtimes_{\alpha}^{\sigma} G$
with the commutant of $A$ in $A \rtimes_{\alpha}^{\sigma} G$ is
nonzero.
\end{thm}

For more details concerning group graded rings in general
and crossed product algebras in particular,
see e.g. \cite{CaenOyst}, \cite{Karp} and \cite{nas}.

Many natural examples of rings, such as rings of
matrices, crossed product algebras defined by separable extensions
and category rings,
are not in any natural way graded by groups, but instead by
categories (see \cite{lu05}, \cite{lu06}, \cite{lu07} and
\cite{oinlun08}).
The main purpose of this article is to obtain a simultaneous
generalization (see Theorem 3) of Theorem \ref{TheoremOne} and Theorem \ref{TheoremTwo}
as well as extending the result from gradings defined by groups
to groupoids.
To be more precise, suppose that $G$ is a category.
The family of objects of $G$ is denoted $\ob(G)$;
we will often identify an object in $G$ with
its associated identity morphism.
The family of morphisms in $G$ is denoted $\ob(G)$;
by abuse of notation, we will often write $s \in G$
when we mean $s \in \mor(G)$.
The domain and codomain of a morphism $s$ in $G$ is denoted
$d(s)$ and $c(s)$ respectively.
We let $G^{(2)}$ denote the collection of composable
pairs of morphisms in $G$ i.e. all $(s,t)$ in
$\mor(G) \times \mor(G)$ satisfying $d(s)=c(t)$.
A category is called a groupoid if all its morphisms are invertible.
Recall from \cite{lu06}
that a ring $R$ is called graded by the category $G$
if there is a set of additive subgroups,
$R_s$, $s \in G$, of $R$ such that $R =
\bigoplus_{s \in G} R_s$ and for all $s ,t \in G$,
we have $R_s R_t \subseteq R_{st}$ if
$(s,t) \in G^{(2)}$ and
$R_s R_t = \{ 0 \}$ otherwise;
if $R_s R_t = R_{st}$,
$(s,t) \in G^{(2)}$, then $R$ is called strongly graded.
By the principal component of $R$ we mean
the set $R_0 := \bigoplus_{e \in \ob(G)} R_e$.
We say that $R$ has the nonzero ideal property
if to each isomorphism $s \in G$
and each nonzero $x \in R_s$, the
right $R_0$-ideal $x R_{s^{-1}}$ is nonzero.
In Section 2, we prove the following result.

\begin{thm}\label{maintheorem}
If a groupoid graded ring has
the nonzero ideal property,
then the intersection of a
nonzero twosided ideal of the ring, with the commutant
of the center of the principal component of the ring, is nonzero.
\end{thm}

Theorem \ref{maintheorem} generalizes Theorem \ref{TheoremOne} and Theorem \ref{TheoremTwo}. In fact, suppose that
$R$ is a ring graded by the group $G$ and that we have chosen $s \in
G$ and a nonzero $x \in R_s$. If $R$ is strongly graded, then $0
\neq x = x 1_R \in x R_{s^{-1}}R_s$, which implies that the right
$R_0$-ideal $x R_{s^{-1}}$ is nonzero. Now suppose that $R$ is a
group graded crossed product $A \rtimes_{\alpha}^{\sigma} G$. Then $x = a_s
u_s$ for some nonzero $a_s \in A_e$. Hence $a_s \alpha(s,s^{-1}) = x
u_{s^{-1}} \in x R_{s^{-1}}$. Therefore the right $R_0$-ideal $x
R_{s^{-1}}$ is nonzero as long as $\alpha(s,s^{-1})$ is not a zero
divisor in $A_e$.

In Section \ref{skewgroupoid}, we generalize \cite[Theorem 3.4]{oinert09}, \cite[Corollary 6]{oinlun08} and \cite[Proposition 10]{oinlun08} by proving the following result.


\begin{thm}\label{skewgroupoidtheorem}
If $A \rtimes^{\sigma} G$ is a skew groupoid algebra with all $A_e$,
for $e \in \ob(G)$, commutative rings and $|\ob(G)|<\infty$, then $A$ is maximal commutative in $A \rtimes^{\sigma} G$
if and only if every intersection of a nonzero twosided ideal of $A \rtimes^{\sigma} G$ and $A$ is nonzero.
\end{thm}

The secondary purpose of this article is to determine the
center of strongly groupoid graded rings $R$
in terms of a groupoid action on $R$
defined by the grading (see Theorem \ref{secondtheorem}
in Section \ref{thecenter}).
This generalizes a result for group graded
rings by the first author together with Silvestrov,
Theohari-Apostolidi and Vavatsoulas (see Lemma 3(iii) in
\cite{oin10}) to the groupoid graded situation.

In Section \ref{examples}, we show that the class of
strongly groupoid graded rings which are not crossed products,
in the sense defined in \cite{oinlun08}, is nonempty for any given groupoid with a finite number of objects. In fact, we show, by an explicit construction,
the following result.

\begin{thm}\label{groupoidexample}
Given a finite groupoid $G$, which has a nonidentity
morphism, there is a ring, strongly graded
by $G$, which is not a crossed product over $G$.
\end{thm}

\section{Ideals}\label{ideals}

In this section, we prove Theorem 3 and a corollary thereof. To this
end, and for use in the next section, we gather some fairly well
known facts from folklore concerning graded rings in a lemma (see
also \cite{lu06} and \cite{nas}). We also show that the commutant of
the principal component of rings graded by cancellable categories,
is again a graded ring (see Proposition (\ref{gradedcommutant}));
this result will be used in Section \ref{skewgroupoid}.

\begin{lem}\label{cancellable}
Suppose that $R$ is a ring graded by a cancellable category $G$.
We use the notation $1_R = \sum_{s \in G} 1_s$
where $1_s \in R_s$, $s \in G$.
(a) $1_R \in R_0$;
(b) if we let $H$
denote the set of $s \in G$ with
$1_{d(s)} \neq 0 \neq 1_{c(s)}$, then $H$
is a subcategory of $G$ with finitely many objects
and $R = \bigoplus_{s \in H} R_s$;
(c) if $G$ is a groupoid (or group), then $H$ is a
groupoid (or group);
(d) if $s \in G$ is an isomorphism, then
$R_s R_{s^{-1}} = R_{c(s)}$
if and only if
$R_s R_t = R_{st}$ for all $t \in G$
with $d(s)=c(t)$.
In particular, if $G$ is a groupoid (or group),
then $R$ is strongly graded if and only if
$R_s R_{s^{-1}} = R_{c(s)}$, $s \in G$.
\end{lem}

\begin{proof}
(a) If $t \in G$, then
$1_t = 1_R 1_t = \sum_{s \in G} 1_s 1_t$.
Since $G$ is cancellable, this implies that
$1_s 1_t = 0$ whenever $s \in G \setminus \ob(G)$.
Therefore, if $s \in G \setminus \ob(G)$, then
$1_s = 1_s 1_R = \sum_{t \in G} 1_s 1_t = 0$.

(b) Since $d(st) = d(t)$, $c(st) = c(s)$
for all $(s,t) \in G^{(2)}$, we get that $H$
is a subcategory of $G$.
By the fact that $1_R = \sum_{s \in ob(H)} 1_s$, we get
that $ob(H)$ is finite. Suppose that $s \in G \setminus H$ is
chosen such that $1_{c(s)} = 0$. Then
$R_s = 1_R R_s = 1_{c(s)} R_s = \{ 0 \}$.
The case when $1_{d(s)} = 0$ is treated similarly.

(c) Suppose that $G$ is a groupoid (or group).
Since $d(s^{-1}) = c(s)$ and $c(s^{-1}) = d(s)$,
$s \in G$, it follows that $H$ is a subgroupoid
(or subgroup) of $G$.

(d) The ''if'' statement is clear. Now we show the
''only if'' statement. Take $(s,t) \in G^{(2)}$ and
suppose that $R_s R_{s^{-1}} = R_{c(s)}$.
Then, by (a), we get that
$R_s R_t \subseteq R_{st} = R_{c(s)} R_{st} =
R_s R_{s^{-1}} R_{st} \subseteq R_s R_{s^{-1}st} =
R_s R_t$. Therefore, $R_s R_t = R_{st}$.
The last part follows immediately.
\end{proof}

\begin{prop}\label{gradedcommutant}
Suppose that $R$ is a ring graded by a category $G$ and that $A$ is
a graded additive subgroup of $R$. For each $s\in G$, denote $C_R(A)_s := C_R(A) \cap R_s$. If $s,t \in G$, then
\begin{itemize}
\item[(a)] $C_R(A)_s  = \bigcap_{u \in G} C_{R_s}(A_u)$;

\item[(b)]
$C_R(A)_s C_R(A)_t \subseteq \left\{
\begin{array}{l}
C_R(A)_{st}, \ {\rm if} \ (s,t) \in G^{(2)}, \\
\{ 0 \}, \ {\rm otherwise;} \\
\end{array}
\right.$

\item[(c)] $C_R(R_0)$ is a
graded subring of $R$ with
$$C_R(R_0)_s =
\left\{
\begin{array}{l}
C_{R_s}(R_{d(s)}), \ {\rm if} \ c(s)=d(s), \\
\left\{ r_s \in R_s \mid R_{c(s)}r_s = r_s R_{d(s)} = \{ 0 \} \right\}, \ {\rm otherwise;} \\
\end{array}
\right.$$

\item[(d)] if $1_R \in R_0$, then $C_R(R_0)$ is a
graded subring of $R$ with
$$C_R(R_0)_s =
\left\{
\begin{array}{l}
C_{R_s}(R_{d(s)}), \ {\rm if} \ c(s)=d(s), \\
\{ 0 \}, \ {\rm otherwise.} \\
\end{array}
\right.$$ In particular, if $G$ is cancellable, then the same
conclusion holds.
\end{itemize}
\end{prop}

\begin{proof}
(a) This is a consequence of the following chain of equalities
$C_R(A)_s = C_R(A) \cap R_s = C_{R_s}(A) = C_{R_s}(\bigoplus_{u \in G}
A_u) = \bigcap_{u \in G} C_{R_s}(A_u)$.

(b) Suppose that $u \in G$, $a_u \in A_u$, $(s,t) \in G^{(2)}$ and
that $r_s \in C_R(A)_s$ and $r_t' \in C_R(A)_t$. Then $(r_s r_t')
a_u = r_s (r_t' a_u) = r_s (a_u r_t') = (r_s a_u) r_t' = (a_u r_s)
r_t' = a_u (r_s r_t')$. Therefore $r_s r_t' \in C_{R_{st}}(A_u)$ for
all $u \in G$. Hence $r_s r_t' \in C_R(A)_{st}$. If, on the other
hand, $(s,t) \notin G^{(2)}$, then, by (a), we get that $C_R(A)_s
C_R(A)_t \subseteq R_s R_t = \{ 0 \}$.

(c) It is clear that $C_R(R_0) \supseteq \bigoplus_{s \in G}
C_R(R_0)_s$. Now we show the reversed inclusion. Take $x \in
C_R(R_0)$, $e \in ob(G)$ and $a_e \in R_e$. Then $\sum_{s \in G} x_s
a_e = \sum_{s \in G} a_e x_s$. By comparing terms of the same
degree, we can conclude that $x_s a_e = a_e x_s$ for all $s \in G$.
Since $e \in ob(G)$ and $a_e \in A_e$ were aritrarily chosen this
implies that $x_s \in C_R(R_0)_s$ for all $s \in G$. Now we show the
second part of (c). Take $e \in ob(G)$. Suppose that $c(s)=d(s)$. If
$d(s) \neq e$, then $C_{R_s}(R_e) = R_s$. Hence $\bigcap_{e \in
ob(G)}C_{R_s}(R_e) = C_{R_s}(R_{d(s)})$. Now suppose that $c(s) \neq
d(s)$. If $c(s) \neq e \neq d(s)$, then $C_{R_s}(R_e) = R_s$.
Therefore, $\bigcap_{e \in ob(G)}C_{R_s}(R_e) = C_{R_s}(R_{c(s)})
\bigcap C_{R_s}(R_{d(s)})$; $C_{R_s}(R_{c(s)})$ equals the set of
$r_s \in R_s$ such that $a r_s = r_s a$ for all $a \in R_{c(s)}$.
Since $d(s) \neq c(s)$, we get that $r_s a_e = 0$;
$C_{R_s}(R_{d(s)})$ is treated similarly.

(d) The claims follow immediately from (c). In fact, suppose that
$c(s) \neq d(s)$. Take $r_s \in R_s$ such that $R_{c(s)}r_s = \{ 0
\}$. Then $r_s = 1_R r_s = 1_{c(s)} r_s = 0$. If $G$ is cancellable,
then, by Lemma \ref{cancellable}(a), the multiplicative identity of
$R$ belongs to $R_0$.
\end{proof}

\subsection*{Proof of Theorem \ref{maintheorem}}
We prove the contrapositive statement. Let $C$ denote the commutant
of $Z(R_0)$ in $R$ and suppose that $I$ is a twosided ideal of $R$
with the property that $I \cap C = \{ 0 \}$. We wish to show that $I
= \{ 0 \}$. Take $x \in I$. If $x \in C$, then by the assumption $x
= 0$. Therefore we now assume that $x  = \sum_{s \in G} x_s  \in I$,
$x_s \in R_s$, $s \in G$, and that $x$ is chosen so that $x \notin
C$ with the set $S := \{ s \in G \mid x_s \neq 0 \}$ of least
possible cardinality $N$. Seeking a contradiction, suppose that $N$
is positive. First note that there is $e \in \ob(G)$ with $1_e x \in
I \setminus C$. In fact, if $1_e x \in C$ for all $e \in \ob(G)$,
then $x = 1_R x = \sum_{e \in \ob(G)} 1_e x \in C$ which is a
contradiction. Note that, by Lemma \ref{cancellable}(b), the sum
$\sum_{e \in \ob(G)} 1_e$, and hence the sum $\sum_{e \in \ob(G)}
1_e x$, is finite. By minimality of $N$, we can assume that
$c(s)=e$, $s \in S$, for some fixed $e \in \ob(G)$. Take $t \in S$.
By the nonzero ideal property there is $y \in R_{t^{-1}}$ with $x_t
y \neq 0$. By minimality of $N$, we can therefore from now on assume
that $e \in S$ and $d(s) = c(s) = e$ for all $s \in S$. Take $d =
\sum_{f \in \ob(G)} d_f \in Z(R_0)$ where $d_f \in R_f$, $f \in
\ob(G)$ and note that $Z(R_0) = \bigoplus_{f \in \ob(G)} Z(R_f)$.
Then $I \ni dx - xd = \sum_{s \in S} \sum_{f \in \ob(G)} (d_f x_s -
x_s d_f)= \sum_{s \in S} d_e x_s - x_s d_e$. In the $R_e$ component
of this sum we have $d_e x_e -x_e d_e=0$ since $d_e \in Z(R_e)$.
Thus, the summand vanishes for $s = e$, and hence, by minimality of
$N$, we get that $dx-xd = 0$. Since $d \in Z(R_0)$ was arbitrarily
chosen, we get that $x \in C$ which is a contradiction. Therefore $N
= 0$ and hence $S = \emptyset$ which in turn implies that $x=0$.
Since $x \in I$ was arbitrarily chosen, we finally get that $I = \{
0 \}$. {\hfill $\square$}

\begin{cor}\label{intersectioncorollary}
If a groupoid graded ring has
the nonzero ideal property and
the principal component of the ring is maximal commutative,
then the intersection of a
nonzero twosided ideal of the ring with
the principal component of the ring is nonzero.
\end{cor}

\begin{proof}
This follows immediately from Theorem \ref{maintheorem}.
\end{proof}

\section{Skew category algebras}\label{skewgroupoid}

We shall recall the definition of a skew category ring from \cite{oinlun08}. By a skew system we mean a
triple $\{ A,G,\sigma \}$ where $G$ is an arbitrary small category, $A$ is the direct sum of rings $A_e$, $e \in \ob(G)$,
$\sigma_s : A_{d(s)} \rightarrow A_{c(s)}$, for $s \in G$,
are ring homomorphisms, satisfying the following two conditions:
\begin{equation}\label{idd}
\sigma_e = {\rm id}_{A_e}
\end{equation}
\begin{equation}\label{algebraa}
\sigma_s(\sigma_t(a)) =  \sigma_{st}(a)
\end{equation}
for all $e \in \ob(G)$, all $(s,t) \in G^{(2)}$ and all $a \in
A_{d(t)}$. Let $A \rtimes^{\sigma} G$ denote the collection of
formal sums $\sum_{s \in G} a_s u_s$, where $a_s \in A_{c(s)}$, $s
\in G$, are chosen so that all but finitely many of them are
zero. Define addition on $A \rtimes^{\sigma} G$ by
(\ref{addition}) and define multiplication on $A \rtimes^{\sigma} G$
as the bilinear extension of the relation
\begin{equation}\label{multiplication}
(a_s u_s)(b_t u_t) = a_s \sigma_s(b_t) u_{st}
\end{equation}
if $(s,t) \in G^{(2)}$ and $(a_s u_s)(b_t u_t) = 0$ otherwise
where $a_s \in A_{c(s)}$ and $b_t \in A_{c(t)}$. One can show that $A \rtimes^{\sigma} G$
has a multiplicative identity if and only if $\ob(G)$
is finite;
in that case the multiplicative identity is
$\sum_{e \in \ob(G)} u_e$.
It is easy to verify that the multiplication on
$A \rtimes^{\sigma} G$ is associative.
Define a left $A$-module structure on
$A \rtimes^{\sigma} G$ by the bilinear extension
of the rule
$a_e (b_s u_s) = (a_e b_s) u_s$
if $e = c(s)$ and
$a_e (b_s u_s) = 0$ otherwise
for all $a_e \in A_e$, $b_s \in A_{c(s)}$,
$e \in \ob(G)$, $s \in G$.
Analogously, define a right $A$-module structure on
$A \rtimes^{\sigma} G$ by the bilinear extension
of the rule
$(b_s u_s) c_f = (b_s \sigma_s(c_f))u_s$
if $f = d(s)$ and
$(b_s u_s) c_f = 0$ otherwise
for all $b_s \in A_{c(s)}$, $c_f \in A_f$,
$f \in \ob(G)$, $s \in G$.
By (\ref{algebraa}) this $A$-bimodule structure
makes $A \rtimes_{\alpha}^{\sigma} G$ an
$A$-algebra. We will often identify $A$ with
$\bigoplus_{e \in \ob(G)} A_e u_e$; this ring will be referred
to as the coefficient ring or principal component of $A \rtimes_{\alpha}^{\sigma} G$.
It is clear that $A \rtimes_{\alpha}^{\sigma} G$ is a category graded ring
in the sense of \cite{lu06} and it is strongly graded.
We call $A \rtimes^{\sigma} G$ the \emph{skew category algebra} associated to the skew system $\{ A,G,\sigma\}$.

\begin{prop}\label{skewcatalg}
If $A \rtimes^{\sigma} G$ is a skew category algebra with all $A_e$,
for $e \in \ob(G)$, commutative rings and $|\ob(G)|<\infty$, such
that every intersection of a nonzero twosided ideal of $A
\rtimes^{\sigma} G$ and $A$ is nonzero, then $A$ is maximal
commutative in $A \rtimes^{\sigma} G$.
\end{prop}

\begin{proof}
We show the contrapositive statement. Suppose that $A$ is not
maximal commutative in $A \rtimes^{\sigma} G$. Then, by Proposition
\ref{gradedcommutant}(d), there exists some $s\in G \setminus G_0$,
with $d(s)=c(s)$, and some nonzero $r_s \in A_{c(s)}$, such that
$r_s u_s$ commutes with all of $A$.
Let $I$ be the (nonzero) ideal in $A \rtimes^{\sigma} G$ generated
by the element $r_s u_{c(s)} - r_s u_s$. Note that all elements of
$I$ are sums of elements of the form
\begin{equation}\label{SpannedBy}
    a_g u_g (r_s - r_s u_s) a_h u_h
\end{equation}
for $g,h\in G$, $a_g\in A_{c(g)}$ and $a_h \in A_{c(h)}$. Suppose
that $(g,h) \in G^{(2)}$, for otherwise the above element is
automatically zero. We may now simplify:
\begin{eqnarray*}
a_g u_g (r_s - r_s u_s) a_h u_h &=&
a_g \sigma_g(r_s a_h) u_{gh} - a_g \sigma_g(a_h r_s) u_{gsh}\\
&=& \underbrace{a_g \sigma_g(r_s a_h)}_{=:b} u_{gh} -
\underbrace{a_g \sigma_g(r_s a_h)}_{=b} u_{gsh}
\end{eqnarray*}
Consider the additive map
\begin{displaymath}
    \varphi : A \rtimes^\sigma G \to A, \quad \sum_{s\in G} a_s u_s \mapsto \sum_{s\in G} a_s .
\end{displaymath}
It is clear that the restriction of $\varphi$ to $A$ is injective.
And since each element of $I$ is a sum of elements of the form
\eqref{SpannedBy} it follows that $\varphi$ is identically zero on
all of $I$. This shows that $I \cap A = \{0\}$ and hence the desired
conclusion follows.
\end{proof}

\subsection*{Proof of Theorem \ref{skewgroupoidtheorem}}
The ''if'' statement follows from Theorem \ref{maintheorem}. The
''only if'' statement follows from Proposition \ref{skewcatalg} if
we let $G$ be a groupoid. {\hfill $\square$}


\section{The center}\label{thecenter}

In this section, we determine the center
of strongly groupoid graded rings (see
Theorem \ref{secondtheorem}) in terms of
an action on the ring induced by the grading
(see Definition \ref{action}).
This is established through
three propositions formulated in a slightly
more general setting.

\begin{prop}\label{firstproposition}
Suppose that $R$ is a ring graded by a category $G$ and that $s \in
G$ is an isomorphism. For each positive integer $i$ take $a_s^{(i)}
, c_s^{(i)} \in R_s$ and $b_{s^{-1}}^{(i)} , d_{s^{-1}}^{(i)} \in
R_{s^{-1}}$ with the property that $a_s^{(i)} = b_{s^{-1}}^{(i)} =
c_s^{(i)} = d_{s^{-1}}^{(i)} = 0$ for all but finitely many $i$. If
$x,y \in C_R(R_{s^{-1}} R_s)$ and
$$A = \sum_{i=1}^{\infty} a_s^{(i)}
xy b_{s^{-1}}^{(i)} \sum_{j=1}^{\infty} c_s^{(j)} d_{s^{-1}}^{(j)}$$
$$B = \sum_{i=1}^{\infty} a_s^{(i)} x b_{s^{-1}}^{(i)}
\sum_{j=1}^{\infty} c_s^{(j)} y d_{s^{-1}}^{(j)}$$ $$C =
\sum_{i=1}^{\infty} a_s^{(i)} b_{s^{-1}}^{(i)} \sum_{j=1}^{\infty}
c_s^{(j)} xy d_{s^{-1}}^{(j)}$$ then $A = B = C$. In particular, if
$G$ is cancellable and
$$\sum_{i=1}^{\infty} a_s^{(i)} b_{s^{-1}}^{(i)} =
\sum_{j=1}^{\infty} c_s^{(j)} d_{s^{-1}}^{(j)} = 1_{c(s)}$$ then the
following equalities hold
\begin{equation}\label{A=C}
\sum_{i=1}^{\infty} a_s^{(i)} x b_{s^{-1}}^{(i)} =
\sum_{j=1}^{\infty} c_s^{(j)} x d_{s^{-1}}^{(j)}
\end{equation}
\begin{equation}\label{A=B}
\sum_{i=1}^{\infty} a_s^{(i)} xy b_{s^{-1}}^{(i)} =
\sum_{i=1}^{\infty} a_s^{(i)} x b_{s^{-1}}^{(i)} \sum_{i=1}^{\infty}
a_s^{(i)} y b_{s^{-1}}^{(i)}
\end{equation}
\end{prop}

\begin{proof}
Suppose that $x,y \in C_R(R_{s^{-1}} R_s)$. The equality $A=B$ (or
$B=C$) follows from the fact that $y$ (or $x$) commutes with
$b_{s^{-1}}^{(i)} c_s^{(j)}$ for all positive integers $i$ and $j$.
The equality (\ref{A=C}) follows from Lemma \ref{cancellable}(a) and the
equality $A=C$ with $y = 1_{d(s)}$. The equality (\ref{A=B}) follows
from Lemma \ref{cancellable}(a), equality (\ref{A=C}) and the
equality $A=B$.
\end{proof}

\begin{prop}\label{secondproposition}
Suppose that $R$ is a ring graded by a category $G$ and that $s,t
\in G$ are isomorphisms with $d(s)=c(t)$. For each positive integer
$i$ take $a_s^{(i)} \in R_s$, $b_{s^{-1}}^{(i)} \in R_{s^{-1}}$,
$c_t^{(i)} \in R_t$, $d_{t^{-1}}^{(i)} \in R_{t^{-1}}$,
$p_{st}^{(i)} \in R_{st}$ and $q_{(st)^{-1}}^{(i)} \in
R_{(st)^{-1}}$ with the property that $a_s^{(i)} = b_{s^{-1}}^{(i)}
= c_t^{(i)} = d_{t^{-1}}^{(i)} = p_{st}^{(i)} = q_{(st)^{-1}}^{(i)}
= 0$ for all but finitely many $i$. If $x \in C_R(R_{(st)^{-1}}
R_s R_t)$ and
$$D = \sum_{k=1}^{\infty} p_{st}^{(k)} x q_{(st)^{-1}}^{(k)}
\sum_{i=1}^{\infty} a_s^{(i)} \sum_{j=1}^{\infty} c_t^{(j)}
d_{t^{-1}}^{(j)} b_{s^{-1}}^{(i)}$$
$$E = \sum_{k=1}^{\infty} p_{st}^{(k)} q_{(st)^{-1}}^{(k)}
\sum_{i=1}^{\infty} a_s^{(i)} \sum_{j=1}^{\infty} c_t^{(j)} x
d_{t^{-1}}^{(j)} b_{s^{-1}}^{(i)}$$ then $D=E$. In particular, if
$G$ is cancellable and the following equalities hold
$$\sum_{k=1}^{\infty} p_{st}^{(k)} q_{(st)^{-1}}^{(k)} =
\sum_{i=1}^{\infty} a_s^{(i)} b_{s^{-1}}^{(i)} = 1_{c(s)} \quad \quad
\sum_{j=1}^{\infty} c_t^{(j)} d_{t^{-1}}^{(j)} =
1_{c(t)}$$ then
\begin{equation}\label{D=E}
\sum_{k=1}^{\infty} p_{st}^{(k)} x q_{(st)^{-1}}^{(k)} =
\sum_{i=1}^{\infty} a_s^{(i)} \sum_{j=1}^{\infty} c_t^{(j)} x
d_{t^{-1}}^{(j)} b_{s^{-1}}^{(i)}
\end{equation}
\end{prop}

\begin{proof}
Suppose that $x \in C_R(R_{(st)^{-1}} R_s R_t)$.
The equality $D=E$ follows from the fact that $x$ commutes with
$q_{(st)^{-1}}^{(k)} a_s^{(i)} c_t^{(j)}$
for all positive integers $i$, $j$ and $k$.
The equality (\ref{D=E}) follows
from Lemma \ref{cancellable}(a) and the equality $D=E$.
\end{proof}

\begin{prop}\label{thirdproposition}
Suppose that $R$ is a ring graded by a category $G$ and that
$s \in G$ is an isomorphism.
For each positive integer $i$ take $a_s^{(i)}
, c_s^{(i)} \in R_s$ and $b_{s^{-1}}^{(i)} , d_{s^{-1}}^{(i)} \in
R_{s^{-1}}$ with the property that $a_s^{(i)} = b_{s^{-1}}^{(i)} =
c_s^{(i)} = d_{s^{-1}}^{(i)} = 0$ for all but finitely many $i$.
If $x \in C_R(R_{s^{-1}} R_{c(s)} R_s)$ and
$y \in R_{c(s)}$, then
\begin{equation}\label{xycommutes}
\sum_{i=1}^{\infty} a_s^{(i)} x b_{s^{-1}}^{(i)}
y \sum_{j=1}^{\infty} c_s^{(j)} d_{s^{-1}}^{(j)} =
\sum_{i=1}^{\infty} a_s^{(i)} b_{s^{-1}}^{(i)} y
\sum_{j=1}^{\infty} c_s^{(j)} x d_{s^{-1}}^{(j)}
\end{equation}
In particular, if
$G$ is cancellable and
$$\sum_{i=1}^{\infty} a_s^{(i)} b_{s^{-1}}^{(i)} =
\sum_{j=1}^{\infty} c_s^{(j)} d_{s^{-1}}^{(j)} = 1_{c(s)}$$ then
\begin{equation}\label{CR}
\sum_{i=1}^{\infty} a_s^{(i)} x b_{s^{-1}}^{(i)}
\in C_R(R_{c(s)})
\end{equation}
If also $x \in Z(R_{d(s)})$, then
\begin{equation}\label{ZR}
\sum_{i=1}^{\infty} a_s^{(i)} x b_{s^{-1}}^{(i)}
\in Z(R_{c(s)})
\end{equation}
\end{prop}

\begin{proof}
Suppose that $x \in C_R(R_{s^{-1}} R_{c(s)} R_s)$ and
$y \in R_{c(s)}$.
The equality (\ref{xycommutes}) follows from the fact that
$x$ commutes with $b_{s^{-1}}^{(i)} y c_s^{(j)} $
for all positive integers $i$ and $j$.
The claim (\ref{CR}) follows from (\ref{A=C}) and
(\ref{xycommutes}).
The claim (\ref{ZR}) follows from (\ref{CR}) and the
fact that $Z(R_e) = R_e \cap C_R(R_e)$ for any
$e \in \ob(G)$.
\end{proof}

\begin{defn}
Suppose that $R$ is a ring strongly graded by a
groupoid $G$.
By abuse of notation, we let $C(R)$
(or $Z(R)$) denote the groupoid with
$C_R(R_e)$ (or $Z(R_e)$), $e \in \ob(G)$, as objects,
and the ring isomorphisms
$C_R(R_{d(s)}) \rightarrow C_R(R_{c(s)})$
(or $Z(R_{d(s)}) \rightarrow Z(R_{c(s)})$),
$s \in G$, as morphisms.
\end{defn}

\begin{defn}\label{action}
Suppose that $R$ is a ring strongly graded by a
groupoid $G$.
For each $s \in G$ and each positive integer $i$, take
$a_s^{(i)} \in R_s$ and
$b_{s^{-1}}^{(i)} \in R_{s^{-1}}$
with the property that
$a_s^{(i)} = b_{s^{-1}}^{(i)} = 0$
for all but finitely many $i$
and $\sum_{i=1}^{\infty} a_s^{(i)} b_{s^{-1}}^{(i)} = 1_{c(s)}$.
Define a function $\sigma_s : R \rightarrow R$ by
$\sigma_s(x) = \sum_{i=1}^{\infty} a_s^{(i)} x b_{s^{-1}}^{(i)}$,
$x \in R$.
By abuse of notation, we let every restriction
of $\sigma_s$ to subsets of $R$ also be denoted by $\sigma_s$.
\end{defn}

\begin{prop}\label{functor}
Suppose that $R$ is a ring strongly graded by a
groupoid $G$. Then the association of each $e \in \ob(G)$
and each $s \in G$ to the ring
$C_R(R_e)$ (or the ring $Z(R_e)$) and the function
$\sigma_s : C_R(R_{d(s)}) \rightarrow C_R(R_{c(s)})$
(or $\sigma_s : Z(R_{d(s)}) \rightarrow Z(R_{c(s)})$)
respectively, defines a functor of groupoids
$\sigma : G \rightarrow C(R)$ (or $\sigma : G \rightarrow Z(R)$).
Moreover, $\sigma$ is uniquely defined on morphisms
by the relations $\sigma_s(x)r_s = r_s x$ and
$\sigma_s(x) 1_{c(s)} = \sigma_s(x)$,
$s \in G$, $x \in C_R(R_{d(s)})$
(or $x \in Z(R_{d(s)})$), $r_s \in R_s$.
\end{prop}

\begin{proof}
We show the claim about $C(R)$.
Since the claim about $Z(R)$ can be shown in a
completely analogous way we leave the details
of this to the reader. Take $s \in G$.
By (\ref{CR}), $\sigma_s$ is well defined.
It is clear that $\sigma_s$ is additive
and that $\sigma_s(1_{R_{d(s)}}) = 1_{R_{c(s)}}$.
By (\ref{A=B}), $\sigma_s$ is multiplicative.
Take $t \in G$ with $d(s)=c(t)$.
By (\ref{D=E}), $\sigma_{st} = \sigma_s \sigma_t$.
By (\ref{A=C}), the definition of $\sigma_s$
does not depend on the particular choice of
$a_s^{(i)} \in R_s$ and $b_{s^{-1}}^{(i)} \in R_{s^{-1}}$
subject to the condition
$\sum_{i=1}^{\infty} a_s^{(i)} b_{s^{-1}}^{(i)} = 1_{c(s)}$.
Therefore, for each $e \in \ob(G)$,
we can make the choice $a_e^{(1)} = b_e^{(1)} = 1_e$
and $a_e^{(i)} = b_e^{(i)} = 0$ for all integers $i \geq 2$;
it is easy to see that this implies that
$\sigma_e = id_{C_R(R_e)}$.
For the second part of the proof suppose that
$s \in G$, $x \in C_R(R_{d(s)})$ and $y \in R$ satisfy
$yr_s = r_s x$ for all $r_s \in R_s$.
Then
$\sigma_s(x) = \sum_{i=1}^{\infty} a_s^{(i)} x b_{s^{-1}}^{(i)}
= \sum_{i=1}^{\infty} y a_s^{(i)} b_{s^{-1}}^{(i)} =
y 1_{c(s)} = y$.
Finally, suppose that $s \in G$, $x \in C_R(R_{d(s)})$
and $r_s \in R_s$. Then
$\sigma_s(x)r_s =
\sum_{i=1}^{\infty} a_s^{(i)} x b_{s^{-1}}^{(i)} r_s =
\sum_{i=1}^{\infty} a_s^{(i)} b_{s^{-1}}^{(i)} r_s x =
1_{c(s)} r_s x = r_s x$. It is clear that
$\sigma_s(x) 1_{c(s)} = \sigma_s(x)$.
\end{proof}

\begin{thm}\label{secondtheorem}
The center of a strongly groupoid graded ring $R$ equals the
collection of $\sum_{e \in \ob(G)} x_e$,
$x_e \in C_R(R_e)$, $e \in \ob(G)$, satisfying
$\sigma_s(x_{d(s)}) = x_{c(s)}$, $s \in G$.
In particular, if $G$ is the disjoint union of
groups $G_e$, $e \in \ob(G)$,
then the center of $R$ equals the direct sum
of the rings $C_R(R_e)^{G_e}$, $e \in \ob(G)$.
\end{thm}

\begin{proof}
Suppose that $y = \sum_{s \in G} y_s$ belongs to the center
of $R$ where $y_s \in R_s$, $s \in G$, and
$y_s = 0$ for all but finitely many $s \in G$.
Since $1_e y = y 1_e$, $e \in \ob(G)$, we get that
$y_s = 0$ whenever $c(s) \neq d(s)$.
Therefore, $y = \sum_{e \in \ob(G)} x_e$ where
$x_e = \sum_{s \in G_e} y_s$, $e \in \ob(G)$.
Since $y \in C_R(R_e)$, $e \in \ob(G)$, we get that
$x_e \in C_R(R_e)$, $e \in \ob(G)$.
Take $s \in G$. The relation
$r_s y = y r_s$, $r_s \in R_s$, and
the last part of Proposition \ref{functor} imply that
$\sigma_s(x_{d(s)}) = x_{c(s)}$.
On the other hand, it is clear, by the
last part of Proposition \ref{functor}, that
all sums $\sum_{e \in \ob(G)} x_e$,
$x_e \in C_R(R_e)$, $e \in \ob(G)$, satisfying
$\sigma_s(x_{d(s)}) = x_{c(s)}$, $s \in G$,
belong to the center of $R$.
The last part of the claim is obvious.
\end{proof}

\section{Examples}\label{examples}

In this section, we show Theorem \ref{groupoidexample}.
Our method will be to generalize, to category graded rings
(see Proposition \ref{categorygradedex}),
the construction given in \cite{das99} for the group graded situation.
To do that, we first need to introduce some more notations.
Let $K$ be a field and $G$ a category. Suppose that $n$ is a
positive integer and choose $s_i \in G$, for $1 \leq i \leq n$. For
$1 \leq i,j \leq n$, let $e_{ij} \in M_n(K)$ be the matrix with 1 in
the $ij$th position and 0 elsewhere. For $s \in G$, we let $R_s$ be
the $K$-vector subspace of $M_n(K)$ spanned by the set of $e_{ij}$,
for $1 \leq i,j \leq n$, such that $(s_i,s) \in G^{(2)}$ and $s_i s
= s_j$.

\begin{prop}\label{categorygradedex}
If $s,t \in G$, then, with the above notations, we get that
\begin{itemize}
\item[(a)] $R_s R_t \subseteq R_{st}$,
if $(s,t) \in G^{(2)}$, and $R_s R_t = \{ 0 \}$, otherwise.

\item[(b)] If $G$ is cancellable, then
the sum $R := \sum_{s \in G} R_s$ is direct. Therefore, in that
case, $R$ is a ring graded by $G$ with components $R_s$, for $s \in
G$.
\end{itemize}
\end{prop}

\begin{proof}
(a) Suppose that $(s,t) \in G^{(2)}$. Take $e_{ij} \in R_s$ and
$e_{lk} \in R_t$. If $j \neq l$, then $e_{ij} e_{lk} = 0 \in
R_{st}$. Now let $j = l$. Then, since $s_i s = s_j$ and $s_j t =
s_k$, we get that $s_i st = s_j t = s_k$. Hence, $e_{ij} e_{jk} =
e_{ik} \in R_{st}$.

(b) Let $X_s$ denote the collection of pairs $(i,j)$, where $1 \leq
i,j \leq n$, such that $(s_i,s) \in G^{(2)}$ and $s_i s = s_j$.
Suppose that $s \neq t$. Seeking a contradiction suppose that $X_s
\cap X_t \neq \emptyset$. Then there would be integers $k$ and $l$,
with $1 \leq k,l \leq n$, such that $s_k s = s_l = s_k t$. By the
cancellability of $G$ this would imply that $s = t$. Therefore, the
sets $X_s$, for $s \in G$, are pairwise disjoint. The claim now
follows from the fact that $R_s = \sum_{(i,j) \in X_s} Ke_{ij}$ for
all $s \in G$.
\end{proof}

\subsection*{Proof of Theorem \ref{groupoidexample}}
Let $H$ be a finite connected groupoid with at least one nonidentity morphism. We begin by showing that one may always find a subring of a matrix ring which is strongly graded by $H$, but which is not a crossed product in the sense of \cite{oinlun08}.
If $H$ only has one object, then it is a group in which case
it has already been treated in \cite{das99}.
Therefore, from now on, we assume that we can choose
two different objects $e$ and $f$ from $H$.
We denote the morphisms of $H$
by $t_1,t_2,\ldots,t_n$. For technical reasons,
we suppose that $t_n = e$ and that $d(t_1) = f$
and $c(t_1)=e$. Let us now choose $n+1$ morphisms $s_1,s_2,\ldots,s_{n+1}$ from $H$
in the following way; $s_i = t_i$, when $1 \leq i \leq n$, and $s_{n+1} = t_n$.

First we define $R$ according to the beginning of Section \ref{examples} and show that it is strongly graded by $H$.
Take $(s,t) \in H^{(2)}$ and $e_{ki} \in R_{st}$.
This means that $s_i st = s_k$. Now pick an integer $j$
with $1 \leq j \leq n$ and $s_i t = s_j$; this is possible
since $\{ s_i \mid 1 \leq i \leq n \} = H$.
Then $e_{ji} \in R_s$ and, since
$s_j s = s_i ts = s_k$, we get that $e_{kj} \in R_s$.
Therefore, $e_{ki} = e_{kj} e_{ji} \in R_s R_t$.
Hence, $R$ is strongly graded.

Now we shall show that $R$ is not a crossed product over $H$
in the sense defined in \cite{oinlun08}.
For $g,h \in \ob(H)$, let $H_{g,h}$ denote the set of
$s \in H$ with $c(s)=g$ and $d(s)=h$.
Since $H$ is connected, all the sets $H_{g,h}$
have the same cardinality; denote this cardinality by $m$.
The component $R_e$ is the $K$-vector space
spanned by the collection of $e_{ij}$ with $s_i e = s_j$,
that is, such that $s_i = s_j$ and $d(s_j)=e$.
Therefore, the $K$-dimension of $R_e$ equals $m+3$.
Furthermore, the component $R_{t_1}$ is the
$K$-vector space spanned by the collection of $e_{ij}$
with $s_i t_1 = s_j$. Since $d(t_1) = f \neq e$,
this implies that the $K$-dimension of $R_{t_1}$
equals $m+1$.
Seeking a contradiction, suppose that $R_{t_1}$ is free on one generator
$u$ as a left $R_e$-module. Then the map $\theta : R_e \rightarrow
R_{t_1}$, defined by $\theta(x) = xu$, for $x \in R_e$, is,
in particular, an isomorphism of $K$-vector spaces.
Since ${\rm dim}_K(R_e) = m+3 > m+1 = {\rm dim}_K(R_{t_1})$,
this is impossible.

We shall now show that our groupoid $G$ is the disjoint union of connected groupoids.
Define an equivalence relation $\sim$ on $\ob(G)$ by saying
that $e \sim f$, for $e,f \in \ob(G)$, if there is a morphism
in $G$ from $e$ to $f$. Choose a set $E$ of representatives
for the different equivalence classes defined by $\sim$.
For each $e \in E$, let $[e]$ denote the equivalence class
to which $e$ belongs. Let $G_{[e]}$ denote the subgroupoid
of $G$ with $[e]$ as set of objects and morphisms $s \in G$
with the property that $c(s),d(s) \in [e]$.
Then each $G_{[e]}$, for $e \in E$, is a connected groupoid
and $G = \biguplus_{e \in E} G_{[e]}$. 
For each $e \in E$, we now wish to define
a strongly $G_{[e]}$-graded ring $R_{[e]}$.
We consider three cases.
If $G_{[e]} = \{ e \}$, then let $R_{[e]} = K$.
If $[e] = \{ e \}$ but the group $G_{[e]}$ contains
a nonidentity morphism, then let $R_{[e]}$
be a strongly $G_{[e]}$-graded ring which is not a crossed
product as defined in \cite{das99}.
If $[e]$ has more than one element, let $R_{[e]}$
denote the strongly $G_{[e]}$-graded ring according to the construction
in the first part of the proof.
We may define a new ring $S$ to be the direct sum $\bigoplus_{e \in E} R_{[e]}$ 
and one concludes that $S$ is strongly graded by $G$ but not a crossed product in the sense of \cite{oinlun08}. 
{\hfill $\square$}

\begin{exmp}
To exemplify Theorem 5, we now give explicitly the construction
in the simplest possible case when $G$ is not a group.
Namely, let $G$ be the unique thin\footnote{In the sense that there is at most one morphism from one object to another.} connected groupoid with two objects.
More concretely this means that the morphisms of $G$ are $e$, $f$,
$s$ and $t$; multiplication is defined by the relations
$$e^2 = e, \quad f^2 = f, \quad es = s, \quad te = t,
\quad sf = s, \quad ft = t.$$
Put $$s_1 = e, \quad s_2 = s, \quad s_3 = t, \quad
s_4 = s_5 = f.$$
and define the $G$-graded ring $R$ as above.
A straightforward calculation shows that
$$
\begin{array}{rcl}
R_e &=& Ke_{11} + Ke_{33} \\
R_f &=& Ke_{22} + Ke_{44} + Ke_{45} + Ke_{54} + Ke_{55} \\
R_s &=& Ke_{12} + Ke_{34} + Ke_{35} \\
R_t &=& Ke_{21} + Ke_{43} + Ke_{53}
\end{array}
$$
Another straightforward calculation shows that
$$R_e R_e = R_e, \quad R_f R_f = R_f, \quad R_e R_s = R_s$$
$$R_t R_e = R_t, \quad R_s R_f = R_s, \quad R_f R_t = R_t.$$
Therefore $R$ is strongly graded by $G$.
However, $R$ is not a
crossed product in the sense defined in \cite{oinlun08}. In fact,
since ${\rm dim}_K(R_f) = 5 > 3 = {\rm dim}_K(R_t)$, the
left $R_f$-module $R_t$ can not be free on one generator.
By a similar argument, the
left $R_e$-module $R_s$ is not free on one generator.
\end{exmp}



%


%



{\bf Acknowledgements:}
The first author was partially supported by The Swedish Research Council, The Crafoord Foundation, The Royal Physiographic Society in Lund, The Swedish Royal Academy of Sciences, The Swedish Foundation of International Cooperation in Research and Higher Education (STINT) and "LieGrits", a Marie Curie Research Training Network funded by the European Community as project MRTN-CT 2003-505078.

\end{document}